\documentclass{amsart}

\usepackage{amsmath}
\usepackage{amssymb}
\usepackage{color}
\usepackage{amsthm}

\theoremstyle{plain}
\newtheorem{theorem}{Theorem}[section]
\theoremstyle{definition}
\newtheorem{definition}{Definition}[section]
\theoremstyle{plain}
\newtheorem{prop}{Proposition}[section]
\theoremstyle{plain}
\newtheorem{lemm}{Lemma}[section]
\theoremstyle{plain}
\newtheorem{cor}{Corollary}[section]

\definecolor{Noir}{rgb}{0,0,0} 
\definecolor{Blanc}{rgb}{1,1,1} 
\definecolor{Gray}{rgb}{0.5,0.5,0.5} 
\definecolor{Rouge}{rgb}{0.8,0.1,0.1} 
\definecolor{DBleu}{RGB}{51,51,178} 
\definecolor{LBleu}{rgb}{0.85,0.85,1} 
\definecolor{Orange}{RGB}{255,140,0}

\newcommand{\bdoc}{\begin{document}} 
\newcommand{\edoc}{\end{document}} 
 
\newcommand{\bcent}{\begin{center}} 
\newcommand{\ecent}{\end{center}} 
 
\newcommand{\benum}{\begin{enumerate}} 
\newcommand{\eenum}{\end{enumerate}} 
\newcommand{\bitem}{\begin{itemize}} 
\newcommand{\eitem}{\end{itemize}} 
 
\newcommand{\btab}{\begin{tabular}} 
\newcommand{\etab}{\end{tabular}} 
\newcommand{\beqn}{\begin{eqnarray}} 
\newcommand{\eeqn}{\end{eqnarray}} 
 
\newcommand{\bmath}{\begin{math}} 
\newcommand{\emath}{\end{math}} 
 
\newcommand{\noin}{\noindent} 
 
\providecommand{\tb}[1]{\textbf{#1}} 
\providecommand{\mb}[1]{\mathbf{#1}} 
 
\newcommand{\bsh}{\backslash} 
 
\newcommand{\ds}{\displaystyle} 
 
\newcommand{\ub}{\underbrace} 
\newcommand{\ob}{\overbrace}

\providecommand{\F}[1]{\mathbb{#1}} 
 
\newcommand{\FF}{\F F} 
\providecommand{\Fn}[1]{\FF_{#1}} 
\newcommand{\Fp}{\Fn{p}} 
\newcommand{\Fq}{\Fn{q}} 
\newcommand{\Fpm}{\Fn{p^m}} 
\newcommand{\Fpn}{\Fn{p^n}} 
\newcommand{\Fpr}{\Fn{p^r}} 
 
\newcommand{\ZZ}{\mathbb{Z}} 
\providecommand{\Zn}[1]{\ZZ_{#1}} 
\providecommand{\ZnZ}[1]{\ZZ/#1\ZZ} 
\newcommand{\Zp}{\Zn{p}} 
\newcommand{\ZpZ}{\ZnZ{p}} 
 
\newcommand{\NN}{\F{N}} 
\newcommand{\QQ}{\F{Q}} 
\newcommand{\RR}{\F{R}} 
\newcommand{\CC}{\F{C}} 
\newcommand{\QQbar}{\overline{\QQ}} 
\newcommand{\Zero}{\mathbb{00}} 
 
\newcommand{\EE}{\F E} 
\newcommand{\II}{\F I} 
\newcommand{\KK}{\F K} 
\newcommand{\MM}{\F M} 
\newcommand{\XX}{\F X} 
\newcommand{\PP}{\F P} 
\newcommand{\FA}{\F A} 
\newcommand{\LL}{\F L} 
\newcommand{\HH}{\mathbb H}  
\newcommand{\FS}{\F S} 
\newcommand{\TT}{\F T}  
\newcommand{\FSig}{\F\Sig} 
\newcommand{\FDel}{\F\Del} 
 
\providecommand{\E}[1]{\hat{\F{#1}}} 
\newcommand{\EC}{\E{C}} 
 
\newcommand{\bQ}{\mathbf{Q}} 
\newcommand{\bP}{\mathbf{P}}

\newcommand{\ind}{\mbox{ind}} 
 
\newcommand{\fx}{f(x)} 
\newcommand{\gx}{g(x)} 
 
\newcommand{\x}{^\star} 
\newcommand{\xs}{^{~\star}} 
\providecommand{\U}[1]{\left(#1\right)\x} 
\newcommand{\xt}{^\times} 
 
\newcommand{\iso}{\simeq} 
\newcommand{\plus}{\oplus} 
\newcommand{\Plus}{\bigoplus} 
\newcommand{\tensor}{\otimes} 
\newcommand{\Tensor}{\bigotimes} 
\newcommand{\inject}{\hookrightarrow} 
\newcommand{\linject}{\hookleftarrow} 
\newcommand{\surject}{\twoheadrightarrow} 
\newcommand{\tri}{\vartriangleleft} 
 
\renewcommand{\ker}{\mbox{ker}\;} 
\newcommand{\Imf}{\mbox{im}~\;} 
\newcommand{\img}{\mbox{im}\;} 
\newcommand{\Hom}{\mbox{Hom}} 
\newcommand{\Sym}{\mbox{Sym}} 
\newcommand{\End}{\mbox{End}\,} 
\newcommand{\Endz}{\mbox{End}_0} 
\newcommand{\Id}{\mbox{Id}} 
\newcommand{\rk}{\mbox{rk}\,} 
\newcommand{\Pic}{\mbox{Pic}} 
\newcommand{\Jac}{\mbox{Jac}} 
\newcommand{\ch}{\mbox{ch}\,} 
\newcommand{\td}{\mbox{td}\,} 
 
\providecommand{\Gal}[1]{\mbox{Gal}(#1)} 
\providecommand{\GAL}[2]{\mbox{Gal}(#1/#2)} 
\providecommand{\Sub}[1]{\mbox{Sub}(#1)} 
\providecommand{\Lat}[1]{\mbox{Lat}(#1)}

\newcommand{\dup}{d_\wedge} 
\newcommand{\drt}{d_>} 

\newcommand{\MF}{\mathfrak}

\newcommand{\Div}{\,\,|\,\,} 
\newcommand{\lcm}{\mbox{lcm}} 
 
\providecommand{\leg}[2]{\left(\frac{#1}{#2}\right)} 
\providecommand{\jac}[2]{\leg{#1}{#2}} 
\providecommand{\qdc}[2]{\left[\frac{#1}{#2}\right]} 
 
\newcommand{\w}{\omega} 
\newcommand{\W}{\Omega} 
\providecommand{\Cal}[1]{\mathcal{#1}} 
\newcommand{\CL}{\Cal L} 
\newcommand{\CO}{\Cal O} 
\renewcommand{\O}{\Cal O} 
\renewcommand{\o}{\Cal O} 
\newcommand{\co}{\Cal O} 
\newcommand{\ca}{\Cal A}
\newcommand{\CR}{\Cal R} 
\newcommand{\CE}{\Cal E} 
\newcommand{\CF}{\Cal F} 
\newcommand{\CQ}{\Cal Q} 
\newcommand{\CK}{\Cal K} 
\newcommand{\CM}{\Cal M} 
\newcommand{\CN}{\Cal N} 
\newcommand{\CDee}{\Cal D} 
\newcommand{\CP}{\Cal P} 
\newcommand{\CS}{\Cal S} 
\newcommand{\ClC}{\Cal C} 
\newcommand{\CJ}{\Cal{J}} 
\newcommand{\CA}{\Cal{A}} 
\newcommand{\CB}{\Cal{B}}

\newcommand{\Si}{\Sigma} 
 
\providecommand{\Ok}[1]{\CO_{#1}} 
\newcommand{\OK}{\Ok{K}}

\renewcommand{\epsilon}{\varepsilon} 
\newcommand{\ep}{\varepsilon} 
 
\providecommand{\abs}[1]{\left|#1\right|} 
\providecommand{\norm}[1]{\lVert#1\rVert} 
 
\newcommand{\di}{\partial} 
\providecommand{\ddy}[1]{\ds\frac{d}{d #1}} 
\providecommand{\didiy}[1]{\ds\frac{\di}{\di #1}} 
\newcommand{\ddx}{\ddy{x}} 
\newcommand{\didix}{\didiy{x}} 
\providecommand{\v}[1]{\vec{#1}} 
 
\newcommand{\ihat}{\hat{\infty}}

\providecommand{\set}[1]{\left\{#1\right\}} 
\providecommand{\lst}[1]{\ds\left[\,\,#1\,\,\right]}

\providecommand{\ip}[2]{\left( #1,#2\right)} 
\providecommand{\IP}[2]{\left\langle #1,#2\right\rangle} 
\providecommand{\bra}[1]{\left\langle\left. #1\right|\right.} 
\providecommand{\ket}[1]{\left.\left| #1\right.\right\rangle} 
\providecommand{\bkv}[4]{\left\langle\left.\tb{#1} #2\right|\tb{#3} #4\right\rangle} 
\providecommand{\bk}[2]{\bkv{}{#1}{}{#2}} 
\providecommand{\bkvop}[5]{\left\langle\tb{#1} #2\left| #5\right|\tb{#3} #4\right\rangle} 
\providecommand{\bkop}[3]{\bkvop{}{#1}{}{#2}{#3}} 
\providecommand{\comm}[2]{\left[#1,#2\right]} 
\providecommand{\expn}[1]{\left\langle #1\right\rangle} 
\providecommand{\gen}[1]{\expn{#1}} 
\newcommand{\Del}{\Delta} 
\newcommand{\Nab}{\nabla} 
\newcommand{\Sig}{\Sigma} 
\newcommand{\oline}{\overline}

\newcommand{\p}{\rho} 
\providecommand{\sprod}[2]{\left\langle #1,#2\right\rangle}

\newcommand{\Ehat}{\overline E} 
\newcommand{\Phihat}{\overline \Phi} 
\newcommand{\phihat}{\overline \phi}

\newcommand\varleq{\mathbin{\vcenter{\hbox{%
  \oalign{\hfil$\scriptstyle<$\hfil\cr 
          \noalign{\kern-.3ex} 
          $\scriptscriptstyle({-})$\cr}%
}}}} 
 
\renewcommand\subsetneq{\mathbin{\vcenter{\hbox{%
  \oalign{\hfil$\scriptstyle\subset$\hfil\cr 
          \noalign{\kern-.3ex} 
          $\scriptscriptstyle({-})$\cr}%
}}}} 

\author{Steven Rayan}
\address{Department of Mathematics \& Statistics\\
  McLean Hall, University of Saskatchewan\\
  Saskatoon, SK, Canada  ~S7N 5E6}
\email{rayan@math.usask.ca}

\title[Quiver at the bottom of the twisted nilpotent cone]{The quiver at the bottom of the twisted nilpotent cone on $\PP^1$}
\date{\today}

\begin{document}

\maketitle

\noin\tb{Abstract.} For the moduli space of Higgs bundles on a Riemann surface of positive genus, critical points of the natural Morse-Bott function lie along the nilpotent cone of the Hitchin fibration and are representations of $\mbox{A}$-type quivers in a twisted category of holomorphic bundles.   The critical points that globally minimize the function are representations of $\mbox{A}_1$.  For twisted Higgs bundles on the projective line, the quiver describing the bottom of the cone is more complicated.  We determine it here.  We show that the moduli space is topologically connected whenever the rank and degree are coprime, thereby verifying conjectural lowest Betti numbers coming from high-energy physics.\\

\section{Introduction}

Let $X$ be a Riemann surface and $\w_X$ its canonical line bundle.  Recall that a \emph{Higgs bundle} on $X$ is a holomorphic vector bundle $E$ adorned with a holomorphic bundle map $\phi:E\to E\tensor\w_X$, usually called a ``Higgs field''.  From the natural K\"ahler metric on the moduli space of stable Higgs bundles over a Riemann surface, we can define a Morse-Bott function $f$ that sends each Higgs bundle $(E,\phi)$ to a scalar multiple of the norm squared of $\phi$.  The existence of this function enables one to use Morse theory to study the topology of the moduli space --- a programme that has been especially successful in low rank (\cite{NJH:86,PBG:94,GGM:07,GHS:14} etc.), but which has been difficult to implement in general.

Emerging from this programme is the now well-known fact that the critical set of $f$ is a submanifold of the \emph{nilpotent cone}, which is precisely the locus of Higgs bundles with nilpotent $\phi$. When the underlying Riemann surface has genus $g\geq1$, $f$ attains an absolute minimum of $0$, and the submanifold of the nilpotent cone along which $f(E,\phi)=0$ is precisely that on which $\phi$ is identically zero.

In other words, the question ``What is the \emph{bottom} of the nilpotent cone?'' has an easy answer when $g\geq1$: it is the moduli space of stable bundles, which is embedded into the moduli space of stable Higgs bundles --- and in particular, into the nilpotent cone --- by $E\mapsto(E,0)$.

In one sense, the same question for $g=0$ has an easy answer, since the moduli space of stable Higgs bundles is empty.  There are two ways to make this less trivial.  One way is to mark the projective line and introduce \emph{parabolic Higgs bundles} adapted to the divisor of the marked points \cite{BY:96}.  This has the advantage that the moduli space is not only nonempty but also hyperk\"ahler, like the moduli space of ordinary Higgs bundles on a positive-genus curve.  The bottom of the cone is now the moduli space of stable parabolic bundles, assuming it is nonempty.

Another way to coax out a nonempty moduli space of stable Higgs bundles on $\PP^1$ is to consider Higgs fields that take values in an arbitrary ample line bundle $\CO(t)$ instead of in $\w$.  The resulting objects are the so-called \emph{twisted Higgs bundles}, studied under various names and mostly in positive genus in \cite{NN:91,EM:94,FB:95,GR:00,GR:16} for example.  At $g=0$, their moduli space lacks the rich hyperk\"ahler structure associated to ordinary and parabolic Higgs bundle moduli spaces. That being said, there is still a natural K\"ahler metric, from which we can define a Morse-Bott function$$f(E,\phi)=\frac{1}{2}\norm{\phi}^2.$$  The moduli space retains a \emph{Hitchin fibration}, which is a proper map to an affine base whose generic fibre is a nonsingular abelian variety.  As in the ordinary Higgs case, the map sends a Higgs bundle to the characteristic polynomial of its Higgs field.  The fibre over zero is an analogue of the global nilpotent variety studied in \cite{VG:01}, containing twisted nilpotent Higgs fields in this case.  However, this twisted nilpotent cone has no natural relationship to the cotangent bundle of a moduli space of bundles.

The answer to the question about the bottom of the cone is not as immediately clear as in the other settings.  First of all, the moduli space of stable bundles is empty when the rank is larger than $1$.  In other words, there are no stable twisted Higgs bundles of the form $(E,0)$, and $f$ does not attain $0$ as its global minimum. (Note that a similar phenomenon occurs for parabolic Higgs bundles with sufficiently small parabolic weights, as in \cite{GGM:07}.)  This leads to a natural question: what are necessary and sufficient conditions on a twisted Higgs bundle $(E,\phi)$ on $\PP^1$ for $\phi$ to minimize $f$?

To answer this, we need not interact directly with $f$.  Rather, we know that global minimizers are exactly the points in the moduli space with \emph{Morse index} $0$, meaning that there are no further downward directions for the Morse flow.  As will be reviewed in $\S\ref{Setup}$, a critical point must have a particular form, called a \emph{holomorphic chain}.  Their appearances in this context the literature include \cite{PBG:95,BGG:04,LAC:09,SSR:11,JHT:12}.  The Morse index can be read off directly from the chain.

In representation-theoretic terms, a holomorphic chain is a representation of an A-type quiver $Q$, but the representations are taken in a category of holomorphic bundles on $X$ with twisted morphisms rather than the category of vector spaces. These objects are a special case of the ``quiver bundles'' considered in \cite{GK:05,LAC:09,LS:12,AS:12}. Here, we consider quivers $Q$ with finite underlying graph $\mbox{A}_n$ for some $n\geq1$. In $Q$, the arrows point in the same direction, from left to right, and we number the nodes sequentially, also from left to right. Before we can specify a representation, we fix an auxiliary bundle $F\to X$ and a labelling of the nodes by pairs of integers $r_i,d_i$ subject to $r_i\geq0$, $\sum r_i=r$, and $\sum d_i=-d$, where $r$ and $d$ are integers for which $0<d<r$:  $$\bullet_{r_1,d_1}\longrightarrow\bullet_{r_2,d_2}\longrightarrow\cdots\longrightarrow\bullet_{r_n,d_n}$$A representation is a $(2n-1)$-tuple$$(U_1,\dots,U_n;\phi_1,\dots,\phi_{n-1})$$ in which $U_i$ is a bundle of rank $r_i$ and degree $d_i$ and $\phi_i$ is an $F$-twisted morphism $\phi_i:U_i\mapsto U_{i+1}\tensor F$.

For the moduli space of ordinary Higgs bundles of rank $r$ and degree $-d$ on a positive-genus curve $X$, we have $F=\w_X$ and the global minimizers of $f$ are representations of the simplest such quiver, $A_1$, with the only possible labelling:$$\bullet_{r,-d}$$ For the correct choice of stability condition, the moduli space of representations associated to this graph is the moduli space of semistable bundles of rank $r$ and degree $-d$. As there are no arrows, the Higgs fields of these bundles are zero.

Finding this quiver for the moduli space of twisted Higgs bundles on $\PP^1$ is equivalent to find its length and labelling.  Fix $F=\CO(t)$ for some $t>0$.  Let $\CM_t(r,-d)$ denote the moduli space of semistable Higgs bundles on $\PP^1$ of rank $r$ and degree $-d$ with $0<d<r$ and Higgs fields taking values in $\CO(t)$; $\mbox{Nilp}_t(r,-d)$, its nilpotent cone; and $\CQ_t(r,-d)$, the moduli space of quiver representations containing the submanifold of $\mbox{Nilp}_t(r,-d)$ along which $f$ is minimized.  We refer to this latter moduli space as a ``quiver-bundle variety'' to avoid confusion with Nakajima quiver varieties.  Our main results are:\\

\noin\tb{Theorem \ref{ThmQuiver}.} When $\gcd(r,d)=1$, $\CQ_t(r,-d)$ is a quiver-bundle variety for the quiver with underlying graph$$\mbox{A}_{\left\lceil\log_{t+1}\left(\frac{r}{d}\right)\right\rceil+1}.$$ If $n=\left\lceil\log_{t+1}\left(\frac{r}{d}\right)\right\rceil+1=2$, its labels are $$\bullet_{r-d,0}\longrightarrow\bullet_{d,-d};$$ if $n=\left\lceil\log_{t+1}\left(\frac{r}{d}\right)\right\rceil+1>2$, then we have$$\bullet_{r-R,0}\longrightarrow\bullet_{dt(t+1)^{n-3},0}\longrightarrow\cdots\longrightarrow\bullet_{dt,0}\longrightarrow\bullet_{d,-d}$$ where $R=d+dt+\cdots+dt(t+1)^{n-3}$. The submanifold of $\CQ_t(r,-d)$ along which $f$ is minimized is the restriction of $\CQ_t(r,-d)$ to the following equivalence classes:\beqn\left\{[(U_1,\dots,U_n;\phi_1,\dots,\phi_{n-1})]\right. & | & U_1,\dots,U_{n-1},U_n\tensor\CO(1)\mbox{ holomorphically trivial,}\nonumber\\ & & \left.\phi_1,\dots,\phi_{n-1}\mbox{ injective}\right\}.\nonumber\eeqn

Above, when we ask for $\phi_i$ to be injective, we mean as a map of global sections.\\

\noin\tb{Theorem \ref{ThmGrTopCon}.}  $\CM_t(r,-d)$ is topologically connected whenever $\gcd(r,d)=1$.\\

For a large range of $r$ and $t$ values that have been inspected by computer, Theorem \ref{ThmGrTopCon} verifies the conjectural lowest Betti numbers for $\CM_t(r,-d)$ coming from Mozgovoy's twisted version of the ADHM recursion formula \cite{SM:12}. These conjectures can presumably be checked using alternative recent results, namely by extracting the Betti numbers from the Donaldson-Thomas invariants for twisted Higgs bundle moduli spaces computed in \cite{MS:14} or by making appropriate modifications to the arguments for ordinary Higgs bundles over finite fields in \cite{OS:16}, so that the closed-form Poincar\'e series obtained in that paper for ordinary Higgs bundles on Riemann surfaces generalizes to $g=0$ and twisted Higgs bundles.  However, it is satisfying to have a direct, Morse-theoretic proof of the connectedness of the twisted Higgs moduli space on $\PP^1$ in the spirit of Hitchin's original approach.\\

\noin\emph{Acknowledgements.} This manuscript was started during the 2016 Symposium on Higgs Bundles in Geometry and Physics held at the Internationales Wissenschaftsforum Heidelberg.  I thank the organizers for their hospitality and for providing a comfortable environment for discussion and work.  I am grateful to Steven Bradlow and Sergey Mozgovoy for useful discussions and comments on the manuscript, and also to Peter Gothen for pointing out a similar phenomenon for small weights in the parabolic case.

\section{Morse theory for twisted Higgs bundles}\label{Setup}

We employ standard notation throughout. In particular, we use $\CO(a)$ to denote a representative of the unique isomorphism class of holomorphic line bundles on $\PP^1$ with degree $a$.  Because the degree map is an isomorphism from the multiplicative group $\mbox{Pic}(\PP^1)$ to the additive group of integers, we have $\CO(a)\tensor\CO(b)\cong\CO(a+b)$.  The dual of a holomorphic vector bundle $E$ is denoted $E^*$.  With these conventions, $\CO(a)^*\cong\CO(-a)$.  By $\End(E)$, we always mean the bundle $E^*\tensor E$.  Its space of global sections, $H^0(\PP^1,\End(E))=H^0(\PP^1,E^*\tensor E)$, is precisely the set of all holomorphic bundle maps from $E$ to itself.  Normally, we omit the $\PP^1$ in sheaf cohomologies $H^i(\PP^1,F)$, as the $\PP^1$ will be understood throughout.

With these conventions, we can formalize what we mean by a twisted Higgs bundle:

\begin{definition} An $\CO(t)$-\emph{twisted Higgs bundle} on $\PP^1$ is a pair $(E,\phi)$ in which $E$ is a holomorphic vector bundle on $\PP^1$ and $\phi$ is an element of $H^0(\End(E)\tensor\CO(t))$.  We refer to the integer $t$ as the \emph{twist} of $(E,\phi)$.\end{definition}

Throughout the paper, $t$ will be a fixed positive integer, and so we can refer to $(E,\phi)$ as a ``twisted Higgs bundle'' without confusion.  It is worth noting that when $t=2$, the Higgs fields are valued in the anticanonical line bundle of $\PP^1$.  These objects are known as \emph{co-Higgs bundles}. They arise in generalized complex geometry and were initially studied in \cite{NJH:10IIa,SSR:11,SSR:13}. 

\begin{definition} A subbundle $U\subset E$ is $\phi$\emph{-invariant} if $\phi(U)\subseteq U\tensor\CO(t)$.  The \emph{slope} of $U$ is the rational number $\mu(U):=\deg(U)/\mbox{rank}(U)$. A twisted Higgs bundle $(E,\phi)$ is called \emph{semistable} if $\mu(U)\leq\mu(E)$ for all nonzero, proper $\phi$-invariant subbundles $U$ of $E$. If the inequality is strict for all such $U$, then $(E,\phi)$ is called \emph{stable}.  If $(E,\phi)$ is not semistable, then it is \emph{unstable}.\end{definition} 

This is Hitchin's slope stability condition from \cite{NJH:86} adapted to the twisted Higgs situation. We denote by $\CM_t(r,d)$ the moduli space of semistable $\CO(t)$-twisted Higgs bundles $(E,\phi)$ on $\PP^1$ in which the rank and degree of $E$ are $r>0$ and $d$, respectively.   It is the set of all semistable pairs $(E,\phi)$ taken up to the following equivalence: $(E,\phi)\cong(E',\phi')$ if there exists a holomorphic bundle isomorphism $\psi:E\to E'$ such that $\phi'=\psi^{-1}\phi\psi$.

We will give a rapid summary of the facts surrounding twisted Higgs bundle moduli spaces on $\PP^1$ and the features of Morse theory that apply to them.
\bitem

\item The set $\CM_t(r,d)$ is nonempty only when $r$ and $t$ are positive. It carries the structure of a smooth, quasiprojective variety of complex dimension $tr^2+1$.  Smoothness is guaranteed by the assumption that $\gcd(r,d)=1$, under which the sets of stable and semistable twisted Higgs bundles of rank $r$ and degree $d$ coincide \cite{NN:91}. The moduli space can be constructed as a GIT quotient \cite{NN:91} or as a K\"ahler quotient and therefore carries a K\"ahler metric (by an adaptation of a construction for ordinary, arbitrary-rank Higgs bundles in \cite{TSCH:88,TSCH:92}).
\item If $d'\cong d\;\mbox{mod}\;r$, then $\CM_t(r,d)$ and $\CM_t(r,d')$ are complex-analytically isomorphic (the map is tensoring by a line bundle of appropriate degree), and so it is enough to work with $\CM_t(r,-d)$ with $d$ in the range $[0,r)$.  Working with nonpositive degree is the convention of our choosing.  As we will eventually restrict to $\gcd(r,d)=1$, we omit $d=0$ and consider $d$ in the open interval $(0,r)$.

\item Consider an $\CO(t)$-twisted Higgs bundle $(E,\phi)$.  The Birkhoff-Grothendieck Theorem, which classifies holomorphic vector bundles on the projective line up to isomorphism, tells us that $E\cong\Plus_{i=1}^r\CO(a_i)$ for some unique set of integers $a_1,\dots,a_r$.  
This also means that $\phi$ is globally represented by an $r\times r$ matrix whose $(i,j)$-th entry is a section $\phi_{ij}$ of the line bundle $\mbox{Hom}(\CO(a_j),\CO(a_i)\tensor\CO(t))=\CO(-a_j+a_i+t)$.  Note that $\phi^*=\overline{\phi^T}$ is a well-defined $\CO(t)$-valued Higgs field for $E^*$, and $(E,\phi)$ is stable if and only if $(E^*,\phi^*)$ is.

\item If a twisted Higgs bundle $(E,\phi)$ is stable, then the kernel of the map $[-,\phi]:H^0(\End(E))\to H^0(\End(E)\otimes\CO(t))$ is $\set{c\mb1_E\,|\,c\in\CC}$.  This is a particular case of the general fact that stable objects are \emph{simple}, meaning that endomorphisms of a stable object are generated by the identity.  In this case, endomorphisms of $(E,\phi)$ are endomorphisms of $E$ that commute with $\phi$.  It follows that, for a stable $(E,\phi)$, the map $[-,\phi]:H^0(\End_0(E))\to H^0(\End_0(E)\otimes\CO(t))$ on trace-free endomorphisms is injective.  By duality, we have that the induced map $[-,\phi]:H^1(\End_0(E))\to H^1(\End_0(E)\otimes\CO(t))$ is surjective.
\item There exists a proper map $h$ from $\CM_t(r,-d)$ to the affine space $B_{t,r}:=\Plus_{i=1}^rH^0(\CO(it))$, sending a twisted Higgs bundle $(E,\phi)$ to the $r$-tuple of characteristic coefficients of $\phi$ (which are sections of various tensor powers of $\CO(t)$).  This map is called the \emph{Hitchin map} or \emph{Hitchin fibration}, first introduced for ordinary Higgs bundles in \cite{NJH:86,NJH:87}.  The fibre $\mbox{Nilp}_t(r,-d):=h^{-1}(0)$ is referred to as the \emph{nilpotent cone}, as it consists of all $(E,\phi)$ for which $\phi$ is nilpotent as a bundle map.

\item The function $f(E,\phi)=\frac{1}{2}\norm{\phi}^2$, where $\norm{\cdot}$ is defined using the K\"ahler metric, is bounded below and is a perfect Morse-Bott function on $\CM_t(r,-d)$.  This fact for twisted Higgs bundles adapts without change from \cite{NJH:86}.
\item The flow of $f$ is coincident with $\mbox{Nilp}_t(r,-d)$ and the critical set of $f$ is a submanifold $\ClC_f\subset\mbox{Nilp}_t(r,-d)$.  Again, this follows without change from properties of ordinary Higgs bundles presented in \cite{TH:98II}.
\item It follows from \cite{PBG:95} that a twisted Higgs bundle $(E,\phi)$ is a critical point of $f$ if and only if $E$ admits a decomposition $E\cong\Plus_{i=1}^n U_i$ for some $n\leq r$, in such a way that $\phi(U_i)\subseteq U_{i+1}\tensor\CO(t)$ for each $i=1,\dots,n-1$ and $\phi(U_n)=0$.  We say that $(E,\phi)$ has the structure of a \emph{holomorphic chain} of \emph{length} $n$.    In particular, $\phi$ is nilpotent of order $n$.  We refer to the bundles $U_i$ as \emph{blocks} of the chain. We say that $\phi$ acts with \emph{weight $1$} on each block.  Note that length $n=1$ is inadmissible when $r>1$, as this would correspond to $(E,0)$, which is always unstable on $\PP^1$.
\item To a holomorphic chain $(E,\phi)$, we assign an $n$-tuple $\mb r=(r_1,\dots,r_n)$ that consists of the ranks of the blocks. 
\item The \emph{Morse index} of $f$ at a critical point $(E,\phi)$ is the number of negative eigenvalues of the Hessian of $f$ at $(E,\phi)$.  Geometrically, this is the number of indepedent downward flow directions of $f$ out of $(E,\phi)$.  Equivalently, if $N_{(E,\phi)}$ is the normal space to $\ClC_f$ at $(E,\phi)$, then the Morse index of $(E,\phi)$ is the dimension of the maximal subspace of $N_{(E,\phi)}$ on which $\mbox{Hess}(f)$ is negative definite.  This number is constant on each connected component of $\ClC_f$.  We will use $\beta(E,\phi)$ to refer to the \emph{complex} Morse index, that is, to the actual Morse index multiplied by $1/2$.  
\item Fix a holomorphic chain $(E,\phi)$ with blocks $U_1,\dots,U_n$, an integer $i$ for which $1\leq i\leq n$, a nonnegative integer $k$, and an integer $q$.  Then, put$$\HH^p_{k,i,q}:=H^p(U_i^*\tensor U_{i+k}\tensor\CO(q))$$ if $k\leq n-i$; otherwise, define it to be the trivial vector space.  We refer to elements of $\Plus_{i=1}^n\HH^p_{k,i,q}$ as $(p,q)$-endomorphisms of $E$ of \emph{weight} $k$.  (In particular, $\phi$ itself is a $(0,t)$-endomorphism of weight $1$.)  Use $h^p_{k,i,q}$ for the complex dimension of this space. 
\item The subspace of $N_{(E,\phi)}$ on which the Hessian of $f$ is negative definite is\beqn\CB(E)\plus\CB(\phi),\nonumber\eeqn where$$\CB(E)=\left\{\begin{array}{ccc}\Plus_{i=1}^{n-1}\Plus_{k=1}^{n-i}\mbox{ker}\left(\HH^1_{k,i,0}\stackrel{[-,\phi]}{\longrightarrow}\HH^1_{k+1,i,t}\right) & \mbox{if} & n>1\\ \set{0} & \mbox{if} & n=0,1\end{array}\right.$$ and $$\CB(\phi)=\left\{\begin{array}{ccc}\Plus_{i=1}^{n-2}\Plus_{k=2}^{n-i}\left(\frac{\HH^0_{k,i,t}}{\mbox{im}\left(\HH^0_{k-1,i,0}\stackrel{[-,\phi]}{\longrightarrow}\HH^0_{k,i,t}\right)}\right) & \mbox{if} & n>2 \\ \set{0} & \mbox{if} & n=0,1,2\end{array}\right.$$ cf. \cite{PBG:95}.  We denote by $\beta(E)$ and $\beta(\phi)$ the complex dimensions of $\CB(E)$ and $\CB(\phi)$, respectively.  In other words, $\beta(E)$ is the subspace of $H^1(\End E)$ consisting of deformations of the complex structure on $E$ that have weight at least $1$ (and for which $\phi$ remains holomorphic), and $\beta(E)$ consists of deformations of $\phi$ of weight at least $2$.  In this language, the direction of the Morse flow is described by weight spaces within the tangent spaces to $\CM_t(r,d)$: the downward flow acts with weight at least $1$ on $E$ and with weight at least $2$ on the Higgs field; the upward flow acts with weight at most $-1$ on $E$ and with weight at most $0$ on the Higgs field.
\item A global minimizer of $f$ is precisely a critical point at which the downward flow terminates.  In other words, a critical point $(E,\phi)$ is a global minimizer of $f$ if and only if $\beta(E,\phi)=0$, which occurs if and only if $\beta(E)=\beta(\phi)=0$.  
\eitem

\noin As we noted earlier, the minimum value of $f$ on $\CM_t(r,-d)$ --- call it $f_{min}$ --- is \emph{positive}, because $\phi$ is never the zero map.  For our purposes, it is easier to classify those $(E,\phi)$ for which $f(E,\phi)=f_{min}$ by using the Morse index. 

\subsection{Calculating the Morse index}\label{SectMorse}

In what follows, let $\delta^m_n$ be $0$ if $m\leq n$ and $1$ otherwise.  Recall that $[-,\phi]:H^1(\End_0(E))\to H^1(\End_0(E)\tensor\CO(t))$ is surjective whenever $(E,\phi)$ is stable.  This forces $$[-,\phi]:\HH^1_{k,i,0}\longrightarrow\HH^1_{k+1,i,t}$$to be surjective, too, whenever $k>0$.  To see this note that, for all $k\geq0$, elements of $\HH^1_{k+1,i,t}=H^1(U_i^*\tensor U_{i+k+1}\tensor\CO(t))$ are trace-free when viewed as twisted endomorphisms of $E$.  Hence, for any $\psi\in\HH^1_{k+1,i,t}$, there exists a $\psi_0\in H^1(\End_0(E))$ for which $[\psi_0,\phi]=\psi$.  It is clear that $\psi_0$ must be an element of $\HH^1_{k,i,0}$, as the action of $[-,\phi]$ always increases weights by exactly $1$ (i.e. a map from $U_i$ to $U_j$ is sent to a map from $U_i$ to $U_{j+1}$), and $[-,\phi]$ simultaneously twists by $\CO(t)$.  

The implication is that$$\beta(E)=\delta^1_n\dim_\CC\Plus_{i=1}^{n-1}\Plus_{k=1}^{n-i}\mbox{ker}\left(\HH^1_{k,i,0}\stackrel{[-,\phi]}{\longrightarrow}\HH^1_{k+1,i,t}\right)=\delta^1_n\sum_{i=1}^{n-1}\sum_{k=1}^{n-i}(h^1_{k,i,0}-h^1_{k+1,i,t}).$$

Each of the differences $h^1_{k,i,0}-h^1_{k+1,i,t}$ is of course nonnegative, as it is the dimension of a subspace of $\HH^1_{k,i,0}$.

For $\beta(E)$, the above expression is all that we will need.  For $\beta(\phi)$, we will need a somewhat finer formula for our arguments.  Consider the unique Birkhoff-Grothendieck decomposition of $E$:$$E\cong\Plus_{1=i}^r\CO(a_i)$$ for some integers $a_i$ such that $\sum a_i=-d$. After re-indexing the integers so that $a_i\geq a_{i+1}$, we denote this non-increasing sequence as $\mbox{BG}(E)$.  If this sequence is equal to $0,\dots,0,-1,\dots,-1$, where the number of $-1$'s in the sequence is $d$, then we say that $E$ has the \emph{generic type}.  (Note that such a sequence is well-defined because $0<d<r$.)

In particular, if $E\cong\Plus_{i=1}^nU_i$, then each summand $U_i$ has its own BG sequence, and the concatenation of $\mbox{BG}(U_1),\dots,\mbox{BG}(U_n)$ is a permutation of $\mbox{BG}(E)$.  If $\mbox{BG}(U_i)=(b_{1,i},\dots,b_{r_i,i})$, then $$\HH^p_{k,i,q}=H^p(U_i^*\tensor U_{i+k}\tensor\CO(q))\cong\Plus_{j=1}^{r_i}H^p(\CO(-b_{j,i})\tensor U_{i+k}\tensor\CO(q)).$$ We will use $b^p_{k,j,i,q}$ for $\dim_\CC H^p(\CO(-b_{j,i})\otimes U_{i+k}\tensor\CO(q))$ (and set this to $0$ when $i+k>n$). Also, if $u\neq j$ or $v\neq i+k-1$, then the image of the map$$[-,\phi]:H^p(\CO(b_{u,i})\tensor U_v\tensor\CO(q)\longrightarrow H^p(\CO(b_{u,i})\tensor U_{v+1}\tensor\CO(q+t))$$ has empty intersection with $H^p(\CO(b_{j,i})\tensor U_{i+k}\tensor\CO(q+t))$.

These considerations imply that$$\Plus_{i=1}^{n-2}\Plus_{k=2}^{n-i}\left(\frac{\HH^0_{k,i,t}}{\mbox{im}\left(\HH^0_{k-1,i,0}\stackrel{[-,\phi]}{\longrightarrow}\HH^0_{k,i,t}\right)}\right)$$ is equal to$$\Plus_{i=1}^{n-2}\Plus_{k=2}^{n-i}\Plus_{j=1}^{r_i}\left(\frac{H^0(\CO(-b_{j,i})\otimes U_{i+k}\otimes\CO(t))}{\mbox{im}\left(H^0(\CO(-b_{j,i})\otimes U_{i+k-1})\stackrel{[-,\phi]}{\longrightarrow}H^0(\CO(-b_{j,i})\otimes U_{i+k}\otimes\CO(t))\right)}\right).$$

Also, note that$$[-,\phi]:H^0(\CO(-b_{j,i})\otimes U_{i+k-1}){\longrightarrow}H^0(\CO(-b_{j,i})\otimes U_{i+k}\otimes\CO(t))$$ is injective since $(E,\phi)$ is stable.

It follows that$$\beta(\phi)=\delta^2_n\sum_{i=1}^{n-2}\sum_{k=2}^{n-i}\sum_{j=1}^{r_i}(b^0_{k,j,i,t}-b^0_{k-1,j,i,0}).$$

Again, note that each term $b^0_{k,j,i,t}-b^0_{k-1,j,i,0}$ is nonnegative, since each is the difference in dimension of two vector spaces, one of which is a subspace of the other.  Hence, $\beta(\phi)=0$ if and only if $$b^0_{k,j,i,t}-b^0_{k-1,j,i,0}=0$$ for all $i,j,k$ in the appropriate ranges. 

\section{Bundle type at the bottom}

From now on, $(E,\phi)$ is a stable holomorphic chain with $\beta(E)=\beta(\phi)=0$; blocks $U_1,\dots,U_n$; and $t>0$ and $0<d<r$.

We also let $u_m$ and $v_m$ stand for the number of $0$'s and $-1$'s, respectively, in the BG sequence of $U_m$. 

\begin{prop}\label{PropBun} $E=\Plus_{i=0}^nU_i$ has the generic type, $U_1\cong\CO^{\plus r_1}$, and $U_n\cong\CO(-1)^{\plus r_n}$, where $r_n\leq d$.\end{prop}

\noin\emph{Remark.} The claim implies that $r_m=u_m+v_m$ for each block and in particular $u_1=r_1,\;v_1=0,\;u_n=0,\;v_n=r_n$.

\begin{proof} Since a chain length of $1$ is prohibited by stability when $r>1$, we can assume that $n\geq2$.  For $(E,\phi)$ to be a global minimizer, we must have $\beta(E)=0$, which is equivalent to$$\sum_{i=0}^{n-1}\sum_{k=1}^{n-i}h^1(U_i^*\tensor U_{i+k})=\sum_{i=0}^{n-1}\sum_{k=1}^{n-i}h^1(U_i^*\tensor U_{i+k+1}\tensor\CO(t)).$$  By Serre duality, the left side of this is$$\CL:=\sum_{i=0}^{n-1}\sum_{k=1}^{n-i}h^0(U_i\tensor U_{i+k}^*\tensor\CO(-2))$$ and the right-hand side is$$\CR:=\sum_{i=0}^{n-1}\sum_{k=1}^{n-i}h^0(U_i\tensor U_{i+k+1}^*\tensor\CO(-t-2)).$$  For purposes of comparison, we re-index the latter sum:$$\CR:=\sum_{i=0}^{n-2}\sum_{k=2}^{n-i}h^0(U_i\tensor U_{i+k}^*\tensor\CO(-t-2)).$$ Note that $\CR$ is nonzero precisely when there are an $i$ and a $k$ such that $U_i$ contains a sub-line bundle isomorphic to $\CO(a)$ and $U_{i+k}$ contains a sub-line bundle isomorphic to $\CO(b)$ with $a\geq b+t+2$.  Since $t$ is positive, this means that $a>b+2$.  In turn, this implies that$$h^0(U_i\tensor U_{i+k}^*\tensor\CO(-2))>h^0(U_i\tensor U_{i+k}^*\tensor\CO(-t-2))>0,$$ and so $\CL$ will not only be positive, but also \emph{strictly larger} than $\CR$.  In other words, $\CL=\CR$ if and only if $\CL=\CR=0$.  Moreoever, since $\CL$ is larger than $\CR$ whenever $\CR$ is nonzero, this means that it is sufficient for us to determine the bundles $E$ for which $\CL=0$.

It follows that a necessary condition for having $\CL=0$ is that for each $e$ in the BG decomposition of $U_i$ and each $f$ in the decomposition of $U_j$ with $j\geq i+1$, we must have $e\leq f+1$.

The block $U_n$ is annihilated by $\phi$ and hence is $\phi$-invariant, and so every sub-line bundle of $U_n$ must have slope less than $-d/r$.  The slope of a sub-line bundle is just its degree, and so we must have $f\leq-1$ for each $f$ in the BG decomposition of $U_n$.  If some $f$ in the decomposition of $U_n$ is less than or equal to $-2$, then this every number in the BG decomposition of $E$ is less than or equal to $-2+1=-1$ (by the condition for having $\CL=0$), which contradicts the fact that $0<d<r$.

Hence, every number in the BG decomposition of $U_n$ is $-1$, which proves the part of the proposition that says that $U_n\cong\CO(-1)^{\plus r_n}$.

Comparing $U_1$ and $U_n$, the $\CL=0$ condition says that every number $e$ in the BG decomposition of $U_1$ must satisfy $e\leq -1+1=0$. In the dual chain $(E^*,\phi^*)$, which has slope $d/r$, $U_1^*$ is the block that is annihilated by $\phi^*$, and so we must have $-e<d/r$ for all $e$ in the decomposition of $U_1$.  If $e<0$, then $-e<d/r$ forces $d/r$ to be larger than $1$, which is a contradiction.  Hence, every number in the BG decomposition of $U_1$ is $0$, which proves the part of the proposition that says that $U_1\cong\CO^{\plus r_1}$ for some $r_1<r$.

If $n=2$, then we are done, as this implies that $r_1+r_2=r$ and the BG decomposition of $E$ is $0,\dots,0,-1,\dots,-1$ (with $r_1$-many $0$'s and $r_2$-many $-1$'s).  If $n>2$, then let $U_j$ be any block with $1<j<n$.  First of all, if $e$ is a number in the BG decomposition of $U_j$, then we must have $e\leq -1+1=0$ (by invoking the $\CL=0$ condition and comparing to $U_n$).  If $e$ were at most $-2$, then every number in the decomposition of $U_1$ would necessarily be bounded above by $-1$, which contradicts the fact that the decomposition of $U_1$ contains only zeroes.  Hence, we must have that $e$ is either $-1$ or $0$.

Hence, every number in the BG decomposition of $E$ is either a $-1$ or a $0$, and so there must be exactly $d$-many $-1$'s.  In other words, $E$ must be of generic type.  This also forces the number $r_n$ to be less than or equal to $d$.

\end{proof}

We sharpen the description of the blocks a little more now.

\begin{lemm}\label{LemmNo-1s}  Let $n\geq4$.  Then $v_1,\dots,v_{n-3}=0$, and $v_{n-2}$ and $v_{n-1}$ cannot be simultaneously nonzero.  Moreover, when $v_{n-2}\neq0$, we must have $t=1$.\end{lemm}

\begin{proof} First of all, by Proposition \ref{PropBun}, we know that $v_1=0$ always.  We begin by assuming $n>4$ and choose any $U_m$ with $m$ in $2\leq m\leq n-3$.  First, because $\beta(\phi)=0$, we must have$$b^0_{j,1,n-1,t}-b^0_{j,1,n-2,0}=0$$ for each $j$ in $1\leq j\leq r_1$, and so$$h^0(\CO(-b_{j,m})\tensor U_n\tensor\CO(t))-h^0(\CO(-b_{j,m}\tensor U_{n-1})=0.$$   Since $U_1\cong\CO^{\plus r_1}$, $U_{n-1}=\CO^{\plus u_{n-1}}\plus\CO(-1)^{\plus v_{n-1}}$, and $U_n\cong\CO(-1)^{\plus r_n}$, the equation becomes$$r_n(-0+(-1)+t+1))-u_{n-1}(-0+0+1)-v_{n-1}(-0+(-1)+1)=0,$$ from which we get $u_{n-1}=r_nt$.  It is also necessary that$$b^0_{j,m,n-m,t}-b^0_{j,m,n-m-1,0}=0$$ for each $j$ in $1\leq j\leq r_m$.  Assume that the BG sequence of $U_m$ contains a $-1$ and choose $j$ so that $b_{j,m}=-1$. The condition$$b^0_{j,m,n-m,t}-b^0_{j,m,n-m-1,0}=0$$ becomes$$r_n(-(-1)+(-1)+t+1))-u_{n-1}(-(-1)+0+1)-v_{n-1}(-(-1)+(-1)+1)=0.$$  Combining this with $u_{n-1}=r_nt$, we obtain $v_{n-1}=r_n(1-t)$.  This forces $v_{n-1}=0$ and $t=1$.

Two additional conditions for $\beta(\phi)=0$ are $$b^0_{j,1,n-2,t}-b^0_{j,1,n-3,0}=0$$ and $$b^0_{j,m,n-m-1,t}-b^0_{j,m,n-m-2,0}=0.$$  Now that $t=1$, $v_{n-1}=0$, and $u_{n-1}=r_n$, the first of these conditions becomes $u_{n-2}=2r_n$ and the second becomes $3u_{n-1}=2u_{n-2}+v_{n-2}$.  Combining them, we get $v_{n-2}=-r_n$, which is a contradiction since $r_n>0$.  Hence, $\mbox{BG}(U_m)$ cannot contain a $-1$ if $2\leq m\leq n-1$.

In the case of $n=4$, we only have $v_{n-1}=v_3$ and $v_{n-2}=v_2$ to be concerned with.  Assume $v_2\neq0$.  The conditions $b^0_{j,1,n-1,t}-b^0_{j,1,n-2,0}=0$ for any $j$ in $1\leq j\leq r_1$ and $b^0_{\ell,2,n-2,t}-b^0_{\ell,2,n-3,0}=0$ for any $\ell$ for which $b_{\ell,j}=-1$, we get $v_3=r_4(1-t)$ which implies that $v_3=0$ and $t=1$.

\end{proof}

An immediate consequence of Lemma \ref{LemmNo-1s} is that either $d=r_n+v_{n-1}$ or $d=r_n+v_{n-2}$ when $r\geq4$, $d=r_3+v_2$ when $r=3$, and $d=r_n$ when $r=2$.

Ruling out the positivity of $v_{n-1}$ and $v_{n-2}$ requires more work and uses stability.

\section{Higgs fields at the bottom}

We are now prepared to prove one of the main theorems.  Note that there exists a unique nonnegative integer $N$ determined by $r$: if $d<r\leq d+dt$, then $N=0$; otherwise, $N$ is the positive integer for which $$d+dt+\cdots+dt(t+1)^{N-1}<r\leq d+dt+\cdots+dt(t+1)^{N}.$$  %
The proof of Theorem \ref{ThmHiggs} uses this number and is divided into cases, but the idea is the same in each case: to locate a destabilizing $\phi$-invariant subbundle whenever $\phi$ does not take a particular form.

\begin{theorem}\label{ThmHiggs} If $d< r\leq d+dt$, then $$\mb r=(r-d,d).$$  If $r>d+dt$, then $$\mb r=(r-R,dt(t+1)^{n-3},\dots,dt(t+1),dt,d),$$ where$$R=d+dt+dt(t+1)+\cdots+dt(t+1)^{n-3}$$ and $$r-R\leq dt(t+1)^{n-2}.$$\end{theorem}

\begin{proof} Recall by Lemma \ref{LemmNo-1s} that $v_{n-1}$ and $v_{n-2}$ are the only numbers $v_m$ that can be nonzero.  We start by assuming that in any case where $v_{n-2}$ is well-defined, i.e. for $n\geq3$, we have $v_{n-2}=0$.  (When $n=3$, $v_{n-2}=v_1=0$ by Proposition \ref{PropBun} directly, but for $n\geq4$ this is not yet obvious.)  This means that $d=r_n+v_{n-1}$.

Now assume that $(E,\phi)$ has $n\geq N+3\geq4$ blocks.   It is a consequence of $\beta(\phi)=0$ that$$u_{n-1}=r_nt,\;u_{n-2}=r_nt(t+1)+v_{n-1}t,\dots,u_{n-(N+1)}=r_nt(t+1)^N+v_{n-1}t(t+1)^{N-1}.$$  The space of global sections of $U_{n-1}\tensor\CO(t)$ can, after a choice of isomorphism, be identified with $\CC^{r_nt(t+1)}\plus\CC^{v_{n-1}t}$.  The space of global sections of $U_{n-2}$ can be identified with the same vector space.  Note that $\phi_{n-2}$, considered as a map of global sections, must be injective; otherwise, its kernel is the space of sections of an invariant, destabilizing trivial subbundle. Hence, $\phi_{n-2}$, as a map of global sections, is an isomorphism of vector spaces.  Now, consider the subbundle $U_{n-1}'$ of $U_{n-1}$ that is isomorphic to $\CO^{\plus r_nt}$, and take $\phi^{-1}(H^0(U_{n-1}'\otimes\CO(t)))$, which is a vector subspace of $H^0(U_{n-2})$ of dimension $r_nt(t+1)$.  Since $U_{n-2}$ contains no positive-degree subbundles, $\phi^{-1}(H^0(U_{n-1}'\otimes\CO(t)))$ must be the space of sections of a subbundle isomorphic to $\CO^{\plus r_nt(t+1)}$. Let $U_{n-2}'$ be this subbundle.  Continue in this way by defining $U_{n-m}'\cong\CO^{\plus r_nt(t+1)^{m-1}}$ to be the subbundle of $U_{n-m}$ whose space of global sections is the preimage of $H^0(U_{n-m+1}'\tensor\CO(t))$ under $\phi_{n-m}$, which again is injective as a map of global sections.  In this way, we get a proper $\phi$-invariant subbundle $U$ of $E$:$$U=U_n\plus U_{n-1}'\plus\cdots\plus U_{n-(N+1)}',$$ which is isomorphic to$$\CO(-1)^{\plus r_n}\plus\CO^{\plus r_nt}\plus\cdots\plus\CO^{\plus r_nt(t+1)^N}.$$  Its slope is$$\mu(U)=\frac{-r_n}{r_n+r_nt+\cdots+r_nt(t+1)^N}=\frac{-d}{d+dt+\cdots+dt(t+1)^N}\geq\frac{-d}{r},$$ and so $U$ is destabilizing.  This bundle is always proper and destabilizing when $n\geq N+3$, even if $v_{n-1}=0$, and so it follows that $n$ can be most $N+2$.  If $n=N+1$, then $u_2=u_{n-(N-1)}=r_nt(t+1)^{N-2}+v_{n-1}t(t+1)^{N-3}$ and by the injectivity of $\phi_1$, $r_n\leq r_nt(t+1)^{N-1}+v_{n-1}t(t+1)^{N-2}$, and so we have a contradiction with $d+dt+\cdots+dt(t+1)^{N-1}<r$.

Hence, $n=N+2$.  Again, $\phi_1$ is necessarily injective and so$$r_1\leq r_nt(t+1)^N+v_{n-1}t(t+1)^{N-1}.$$  Recalling that $d+dt+\cdots+dt(t+1)^{N-1}<r$, we must also have $r_1>v_{n-1}t(t+1)^{N-1}$.  However, because $r\leq d+dt+\cdots+dt(t+1)^N$, we have a contradiction if $r>d+dt+\cdots+dt(t+1)^{N-1}+r_nt(t+1)^N$, which is resolved only if $v_{n-1}=0$.  So now restrict to the range $r\leq d+dt+\cdots+d(t+1)^{N-1}+r_nt(t+1)^N$. We can write $r_1$ as $K+v_{n-1}t(t+1)^{N-1}$ where $1\leq K\leq r_nt(t+1)^N$.  Define in the same way as above a subbundle $U$, but take $U_1'\cong\CO^{\plus K}$ to be a subbundle of $U_1$ whose global sections lie in the preimage of $H^0(U_2'\tensor\CO(t))\cong\CC^{r_nt(t+1)^N}$ under $\phi_1$.  Note that $K/r_n>K/d$ and that$$\mu(E)=\frac{-d}{r}=\frac{-d}{d+dt+\cdots+dt(t+1)^{N-1}+K}=\frac{-1}{1+t+\cdots+t(t+1)^{N-1}+\frac{K}{d}}.$$  Comparatively,$$\mu(U)=\frac{-r_n}{r_n+r_nt+\cdots+r_nt(t+1)^{N-1}+K}=\frac{-1}{1+t+\cdots+t(t+1)^{N-1}+\frac{K}{r_n}}.$$  It follows immediately that $U$ is destabilizing, unless $v_{n-1}=0$, in which case $U$ is no longer proper.

Now we have $v_{n-1}=0$ and it now follows that $r_m=u_m$ for $1\leq m\leq n-1$ and $r_n=d$, and so $\mb r$ is as described in the statement of the theorem.

Finally, we want to eliminate completely the case of $v_{n-2}>0$ when $n\geq4$.  To do this, we assume $n\geq4$ and that $v_{n-2}$ is nonzero and seek contradictions.  First, we must have $v_{n-1}=0$ and $t=1$ by Lemma \ref{LemmNo-1s}, and so $r$ lies in either the range$$d<r\leq2d,$$ for which $N=0$, or$$d+2^0d+\cdots+2^{N-1}d<r\leq d+2^0d+\cdots+2^Nd.$$ Assume that $n\geq N+3$. The $\beta(\phi)=0$ conditions manifest themselves as$$u_{n-1}=r_n,\;u_{n-2}=2r_n,\;u_{n-3}=4r_n+v_{n-2},\dots,u_{n-(N+1)}=2^Nr_n+2^{N-2}v_{n-2}.$$  We define a proper $\phi$-invariant subbundle$$U=U_n\plus U_{n-1}\plus U_{n-2}'\plus\cdots\plus U_{n-(N+1)}',$$ where $U_{n-2}'$ is the subbundle of $U_{n-2}$ isomorphic to $\CO^{\plus 2r_n}$, and further $U_{n-m}'$ are defined as above to be the subbundle of $U_{n-m}$ whose space of global sections is the preimage of those of $U_{n-m+1}'\tensor\CO(t)$ under $\phi_{n-m}$.  As above, $U_{n-m}'\cong\CO^{\plus 2^{m-1}r_n}$.  The slope is$$\mu(U)=\frac{-r_n}{r_n+r_n+2r_n+\cdots+2^Nr_n}=\frac{-d}{d+2^0d+\cdots+2^Nd}\geq\frac{-d}{r}$$ and so this bundle is destabilizing.  This bundle is always proper and destabilizing when $n\geq N+3$, even if $v_{n-2}=0$, and so we must have $n\leq N+2$. As in the preceding arguments, it is easy to eliminate values of $n$ smaller than $N+2$.  Since $n\geq4$ by assumption, we must have $N$ at least $2$.

It follows from the injectivity of $\phi_1$ that $r_1\leq 2^Nr_n+2^{N-2}v_{n-2}$.  This means that we cannot have$$r>d+2^0d+\cdots+2^{N-2}d+2^{N-1}r_n+2^Nr_n,$$and so for such $r$, we have a contradiction.  So now we restrict to$$d+2^0d+\cdots+2^{N-1}d<r\leq d+2^0d+\cdots+2^{N-2}d+2^{N-1}r_n+2^Nr_n.$$This provides the further restriction $2^{N-1}d<(2^{N-1}+2^N)r_n$, which is equivalent to $d<3r_n$. We can write $r_1=K+2^{N-2}v_{n-2}$ where $1\leq K\leq2^Nr_n$.  Define in the same way as above a subbundle $U$, but take $U_1'\cong\CO^{\plus K}$ to be a subbundle of $U_1$ whose global sections lie in the preimage of $H^0(U_2'\tensor\CO(1))\cong\CC^{2^Nr_n}$ under $\phi_1$.  Note that $K/r_n>(K-2^{N-1}v_{n-2})/d$ and that\beqn\mu(E) & = & \frac{-d}{d+2^0d+\cdots+2^{N-1}d+(K-2^{N-1}v_{n-2})}\nonumber\\ & = & \frac{-1}{1+t+\cdots+t(t+1)^{N-1}+\frac{K}{d}}.\nonumber\eeqn  On the other hand,$$\mu(U)=\frac{-r_n}{r_n+r_nt+\cdots+r_nt(t+1)^{N-1}+K}=\frac{-1}{1+t+\cdots+t(t+1)^{N-1}+\frac{K}{r_n}}.$$  It follows immediately that $U$ is destabilizing, unless $v_{n-2}=0$, in which case $U$ is no longer proper.

Having eliminated $v_{n-2}>0$ in every possible case, we default to the preceding arguments and so $\phi$ has the claimed shape.

\end{proof}

\noin\emph{Remark.} The $d=1$ case of this result appears in \cite{SSR:11}.  The theorem for that particular case is markedly easier to establish than for general $d$.  One reason is that for $d=1$ it is possible to use an inductive argument based on the rank.  If we delete a line bundle from $U_1$ and then restrict the Higgs field to the resulting rank $r-1$ bundle, the new Higgs bundle is a stable minimizer.  Its stability is ensured because the distribution of $-1$'s amongst the blocks is already known: there is only one $-1$ and so it must be in $U_n$ (which is then just $\CO(-1)$ itself) and so every proper invariant subbundle will have slope $-1/r'<-1/(r-1)$.  It is also true for general $d$ coprime to $r$ that this restriction procedure produces successive stable minimizers, but this is only known \emph{a posteriori}, after proving Theorem \ref{ThmHiggs}.\\

With Theorem \ref{ThmHiggs} come the following immediate corollaries:\\

\begin{cor}\label{CorTrivial} If $(E,\phi)$ is a stable global minimizer of $f$, then $U_1,\dots,U_{n-1}$ and $U_n\tensor\CO(1)$ are holomorphically trivial.\end{cor}

\begin{cor}\label{CorUnique} Any two stable global minimizers of $f$ for the same $r$, $d$, and $t$ have equal lengths and equal rank vectors and their corresponding blocks are isomorphic as holomorphic bundles.\end{cor}

Conversely, it is easy to check that any stable $(E,\phi)$ with structure specified by Theorem \ref{ThmHiggs} and Corollary \ref{CorTrivial} satisfies $\beta(E)=\beta(\phi)=0$.  Such a Higgs bundle is stable precisely when the maps $\phi_1,\dots,\phi_{n-1}$ are injective as maps of global sections.  That this is necessary is a consequence of the proof of Theorem \ref{ThmHiggs}.  (There, we saw that $\phi_i$ restricted to any holomorphically-trivial subbundle of $U_i$ must induce an injective map of global sections.  We now know that $U_1,\dots,U_{n-1}$ are holomorphically trivial themselves.)  One can show this is sufficient as well, by using the injectivity in combination with an appropriate automorphism of $U_n$ (applied to $\phi_{n-1}$) to show that every invariant subbundle has maximally-negative degree $-d$.

Although stability has been assumed all along, we include that hypothesis explicitly in Corollary \ref{CorUnique} to emphasize that the statement is not necessarily true when $\gcd(r,d)\neq1$, as there may be minimizers that are semistable but not stable and which do not take the shape prescribed by Theorem \ref{ThmHiggs}. Case in point, consider $r=4$, $d=2$, and $t=1$.  Any critical point $(E,\phi)$ with $U_1\cong\CO$, $U_2\cong\CO\plus\CO(-1)$, and $U_3\cong\CO(-1)$ is at best semistable, as there is always an invariant subbundle of degree $-1/2$.  (This Higgs bundle has $v_{n-1}\neq0$ but does not violate the proof of Theorem \ref{ThmHiggs} because it is not strictly stable.)  An example Higgs field for this bundle that attains semistability is the one that maps $U_1$ identically onto $\CO(-1)\tensor\CO(1)\subset U_2\CO(1)$, $\CO\subset U_2$ identically onto $\CO(-1)\tensor\CO(1)\subset U_3\tensor\CO(1)$, and $\CO(-1)\subset U_2$ to $\CO(-1)\tensor\CO(1)\subset U_3\tensor\CO(1)$ via a choice of nonzero section $p\in H^0(\CO(1))$. One can check that this example satisfies $\beta(E)=0$ and $\beta(\phi)=0$ by appealing directly to their definitions (as opposed to the formulas derived in Section \eqref{SectMorse}, which depend on strict stability).

Returning to strictly stable minimizers, observe that we can express the length $n$ in terms of $r$, $d$, and $t$.  If $N\geq1$, we can take the sum of the geometric series with common ratio $t+1$ and write$$(t+1)^N<\frac{r}{d}\leq(t+1)^{N+1}.$$  Using the fact that $n=N+2$, we get$$ n<\log_{t+1}\left(\frac{r}{d}\right)+2\leq n+1,$$ from which we obtain a closed-form formula for $n$ as a function of $r,d,t$:$$n(r,d,t)=\left\lceil\log_{t+1}\left(\frac{r}{d}\right)\right\rceil+1.$$

Notice that this formula produces all the relevant values of $n$, for when $d<r\leq d+dt$ this formula returns $2$.

Taken all together, we arrive at the following theorem:

\begin{theorem}\label{ThmComplete} Let $t,d,r$ be positive integers for which $0<d<r$. If$$\left\lceil\log_{t+1}\left(\frac{r}{d}\right)\right\rceil+1=2,$$then a stable Higgs bundle $(E,\phi)$ is a global minimizer of $f$ in $\CM_t(r,-d)$ if and only if\bitem\item $E\cong U_1\oplus U_2$\item $\mbox{\emph{rk}}(U_1)=r-d$ and $\mbox{\emph{rk}}(U_2)=d$ \item $U_1$ and $U_2\tensor\CO(1)$ are holomorphically trivial\item $\phi(U_1)\subseteq U_{2}\tensor\CO(t)$ and $\phi$ is injective as a map of global sections.\eitem If$$\left\lceil\log_{t+1}\left(\frac{r}{d}\right)\right\rceil+1=n>2,$$then a stable Higgs bundle $(E,\phi)$ is a global minimizer of $f$ in $\CM_t(r,-d)$ if and only if \bitem\item $E\cong\Plus_{i=1}^nU_i$\item$\mbox{\emph{rk}}(U_n)=d$, $\mbox{\emph{rk}}(U_i)=dt(t+1)^{n-i-1}$ for $2\leq i\leq n-1$, and $\mbox{\emph{rk}}(U_1)=r-(d+dt+dt(t+1)+\cdots+dt(t+1)^{n-3})$\item $U_1\dots,U_{n-1}$ and $U_n\tensor\CO(1)$ are holomorphically trivial\item $\phi(U_i)\subseteq U_{i+1}\tensor\CO(t)$ for $1\leq i\leq n-1$ and $\phi(U_n)=0$\item $\phi_1,\dots,\phi_{n-1}$ are injective as maps of global sections.  If $\gcd(r,d)=1$, then all minimizers are of the form described above.\eitem\end{theorem}

\section{The quiver}

Theorem \ref{ThmComplete} encodes the quiver $\CQ_t(r,-d)$ whose moduli space of representations contains all the stable minimizers of $f$ in $\mbox{Nilp}_t(r,-d)\subset\CM_t(r,-d)$, and also tells us the locus in $\CQ_t(r,-d)$ along which $f(E,\phi)=f_{min}$.  The following is a rephrasing of Theorem \ref{ThmComplete}:

\begin{theorem}\label{ThmQuiver} When $\gcd(r,d)=1$, $\CQ_t(r,-d)$ is a quiver-bundle variety for the quiver with underlying graph$$\mbox{A}_{\left\lceil\log_{t+1}\left(\frac{r}{d}\right)\right\rceil+1}.$$ If $n=\left\lceil\log_{t+1}\left(\frac{r}{d}\right)\right\rceil+1=2$, its labels are $$\bullet_{r-d,0}\longrightarrow\bullet_{d,-d};$$ if $n=\left\lceil\log_{t+1}\left(\frac{r}{d}\right)\right\rceil+1>2$, then we have$$\bullet_{r-R,0}\longrightarrow\bullet_{dt(t+1)^{n-3},0}\longrightarrow\cdots\longrightarrow\bullet_{dt,0}\longrightarrow\bullet_{d,-d}$$ where $R=d+dt+\cdots+dt(t+1)^{n-3}$. The submanifold of $\CQ_t(r,-d)$ along which $f$ is minimized is the restriction of $\CQ_t(r,-d)$ to the following equivalence classes:\beqn\left\{[(U_1,\dots,U_n;\phi_1,\dots,\phi_{n-1})]\right. & | & U_1,\dots,U_{n-1},U_n\tensor\CO(1)\mbox{ holomorphically trivial,}\nonumber\\ & & \left.\phi_1,\dots,\phi_{n-1}\mbox{ injective}\right\}.\nonumber\eeqn\end{theorem}

When $\gcd(r,d)\neq1$ (for example, when $r=d+dt+\cdots+dt(t+1)^N$ and $d>1$), there may be semistable but not stable minimizers that do not correspond to the quiver above.  In the previous section, we saw that $\CM_1(4,2)$ has a minimizer corresponding to the quiver

$$\bullet_{1,0}\longrightarrow\bullet_{2,-1}\longrightarrow\bullet_{1,-1}$$

The key difference is that this quiver has an interior node of nonzero degree.

\section{Components and the ADHM recursion formula}

By Hitchin's Morse-theoretic localization procedure \cite{NJH:86}, the Poincar\'e polynomial for the ordinary rational cohomology of the moduli space is$$\CP(\CM_t(r,-d);y)=\sum_iy^{2\beta(\CN_i)}\CP_y(\CN_i),$$ where $\CN_i$ is the $i$-th connected component of the critical set of $f$, $\beta(\CN_i)$ is the complex Morse index of any Higgs bundle in $\CN_i$, and $\CP_y(\CN_i)$ is the Poincar\'e polynomial of $\CN_i$ in the variable $y$.  It is worth noting that, because of the properness of the Hitchin fibration and hence the compactness of $\mbox{Nilp}_t(r,-d)$, there are only finitely-many critical components.

The constant term of $\CP(\CM_t(r,-d);y)$ is the number of connected components of $\CM_t(r,-d)$.  It is clear that this number of connected components of the set of global minimizers of the Morse-Bott function.  Having classified the solutions of $f(E,\phi)=f_{min}$ when $\gcd(r,d)=1$, we can count these components.  

\begin{theorem}\label{ThmGrTopCon}  Assume $\gcd(r,d)=1$.  As a complex variety, the space of global minimizers of $f$ is isomorphic to the Grassmannian$$\mbox{\emph{Gr}}(r-d,dt)$$when $n=\left\lceil\log_{t+1}\left(\frac{r}{d}\right)\right\rceil=2$, and is isomorphic to$$\mbox{\emph{Gr}}(r-R,dt(t+1)^{n-2}),$$ where $R=d+dt+dt(t+1)+\cdots+dt(t+1)^{n-3}$, when $n=\left\lceil\log_{t+1}\left(\frac{r}{d}\right)\right\rceil+1>2$.  As such, $\CM_t(r,-d)$ is topologically connected.\end{theorem}

\begin{proof} By Theorem \ref{ThmComplete}, the blocks $U_1,\dots,U_n$ of the Higgs bundles that minimize $f$ have fixed holomorphic type, and so the moduli are concentrated in the Higgs fields.  The condition on the Higgs fields is that they are injective as maps of sections.  For $\phi_i:U_i\to U_{i+1}\tensor\CO(t)$ with $i>1$, this requires that $\phi_i$ is an isomorphism of spaces of global sections, as $H^0(U_i)$ and $H^0(U_{i+1}\tensor\CO(t))$ always have the same dimension (also by Theorem \ref{ThmComplete}).  Noting that $\phi_i$ as a bundle map is determined by its induced map on global sections of $U_i\cong\CO^{\oplus r_i}$, we have a one-to-one correspondence between admissible maps $\phi_i$ and elements of $\mbox{GL}(r_i,\CC)$ for each $i>1$.  Each of these maps of global sections is acted on the right by $\mbox{Aut}(U_i)\cong\mbox{GL}(r_i,\CC)$.  The quotient is just $\set{\mbox{I}}\in\mbox{GL}(r_i,\CC)$.  What remains is $\phi_1$.  This is an injective map from $H^0(U_1)$ to $H^0(U_2\tensor\CO(t))$, but now $h^0(U_1)\leq h^0(U_2\tensor\CO(t))$ by Theorem \ref{ThmComplete}.  In other words, $\phi_1$ is in correspondence with embeddings of an $r_1$-plane into $H^0(U_2\tensor\CO(t))$.  Quotienting by the right multiplication action of $\mbox{Aut}(U_1)\cong\mbox{GL}(r_1,\CC)$ gives us a Grassmannian.  This is one of the two Grassmannians in the statement of the theorem, depending on whether $n=2$ (in which case $U_2$ is $\CO(-1)^{\oplus d}$) or $n>2$ (in which case $U_2$ is trivial).\end{proof}

In particular, the submanifold of stable minimizers has dimension$$r_1(dt(t+1)^{n-2}-r_1),$$where $r_1=r-(d+dt+\cdots+dt(t+1)^{n-3})$ when $n>2$, $r_1=r-d$ when $n=2$, and $1\leq r_1\leq dt(t+1)^{n-2}$.  There is a unique stable global minimizer up to isomorphism when$$r=d+dt+\cdots+dt(t+1)^{n-2}.$$ In particular, when $d=1$, there is a unique global minimizer up to isomorphism for each of these ranks.  When $d\neq0$, these ranks are not coprime to $d$, and there will in general be semistable minimizers in addition to the unique stable one.

Finally, we remark that our calculation of the lowest nonzero Betti number of $\CM_t(r,-d)$, which is always $1$ whenever $\gcd(r,d)=1$, verifies the conjectural lowest Betti number coming from the twisted ADHM recursion formula for a large range of ranks and twists that can be checked by computer.  The twisted ADHM recursion formula was posed and studied by Mozgovoy \cite{SM:12} as a generalization of the Chuang-Diaconescu-Pan ADHM recursion formula coming from physics \cite{CDP:11} (see also \cite{CDDP:15}).  Regarding its solutions, we should add that the ADHM Betti numbers depend on two parameters which can be identified with $r$ and $t$.  There is no dependence on $d$, which is consistent with the fact that the Betti numbers of ordinary Higgs bundle moduli spaces are independent of the degree, as proved in \cite{HRV:08}, at least when $\gcd(r,d)=1$.  The proof in \cite{HRV:08} relies on the homeomorphism to a character variety induced by nonabelian Hodge theory, which does not apply for non-parabolic twisted Higgs bundles on $\PP^1$, and so this invariance is merely conjectural for such twisted Higgs bundles.

We should note, however, that we already have an example of an instance where the moduli space is not connected when $r$ and $d$ are not coprime.  As above demonstrated above, $\CM_2(4,2)$ not only has minimizers of the form prescribed by Theorem \ref{ThmComplete}, but also semistable minimizers as elicited in the example following Corollary \ref{CorUnique}.  It is known that the degree of a block $U_i$ is constant on a connected component of the critical set of $f$ (Lemma 9.2 in \cite{HT:03}, attributed to C. Simpson), and so these two types of minimizers cannot belong to the same component, as the degree of $U_2$ and $U_3$ in one do not respectively match those of $U_2$ and $U_3$ in the other.  This disconnects the moduli space, and gives an explicit example of how the degree invariance of the Poincar\'e polynomial fails when $\gcd(r,d)$ is left unrestricted.

\bibliographystyle{acm} 
\bibliography{biblio}

\end{document}